\documentclass{article}

\usepackage{amsmath,amsfonts}

\usepackage{geometry}
\usepackage{amssymb}

\newtheorem{theorem}{Theorem}
\newtheorem{lemma}[theorem]{Lemma}

\newtheorem{remark}{Remark}

\newtheorem{proposition}{Proposition}

\newcommand{\sectionnew}[1]{
\section{#1}\setcounter{equation}{0}
\setcounter{theorem}{0}}

\newcommand{\R}{\mathbb{R}}
\newcommand{\T}{\mathbb{T}}
\newcommand{\ep}{\varepsilon}
\newcommand{\diver}{\rm{div}}
\newcommand{\Tr}{\rm{Tr}}
\newcommand{\para}{{_{||}}}

\newcommand{\cqfd}
{
\mbox{} \nolinebreak  \hfill  \rule{2mm}{2mm} \\
\newline
}

\title{DERIVATION OF A GYROKINETIC MODEL.\\
EXISTENCE AND UNIQUENESS OF SPECIFIC STATIONARY SOLUTIONS.}
\author{PHILIPPE GHENDRIH\\ \\CEA, DRFC\\  F-13108 Saint-Paul lez Durance, France\\
\\
\hspace*{0.05in}MAXIME HAURAY\\
\\LJLL, Paris 6, 7 \& CNRS\\175 rue du Chevaleret, 75013 Paris, France\\
\\
\hspace*{0.05in}ANNE NOURI\\
\\LATP, Universit\'e de Provence\\ 39 rue F.Joliot Curie, 13453
Marseille Cedex 13, France}

\date{}

\begin{document}

\maketitle

{\noindent \bf Abstract.}\hspace{0.1in}\\
A finite Larmor radius approximation is derived from the classical Vlasov equation, in the limit of large (and uniform) external magnetic field. We also provide an heuristic derivation of the  electroneutrality equation in the finite Larmor radius setting. Existence and uniqueness of a solution is proven in the stationary frame for solutions depending only on the direction parallel to the magnetic field and factorizing in the velocity variables.\\

\footnotetext[1]{2000 Mathematics Subject Classification. 41A60, 76P05, 82A70,
78A35.}
\footnotetext[2]{Key words. Plasmas, gyrokinetic model, steady state solutions.}

{\noindent \bf Introduction.}\hspace{0.1in}
The ITER project is a challenge to the growing need of new sources of energy. It aims at producing energy by nuclear fusion. Nuclear reactions take place in a tokamak, where a high temperature plasma is confined. So far, confined plasmas are performed with relatively short energy confinement times due to microscale instabilities that generate turbulent transport \cite{ITER_1999}. It is observed that the characteristic frequencies of these instabilities is several orders of magnitude smaller than the ion Larmor gyration frequency governed by the strong magnetic field. Studies of nuclear fusion in tokamaks are in full expansion, both experimentally and theoretically. Kinetic models are appropriate for studying the core of the plasma since the collisions have a very weak effect in these hot and low density plasmas. Physicists use gyrokinetic models and especially the finite Larmor radius approximation to model the core of the plasma \cite{Garbet_2009}. Taking into account the fast Larmor gyration of the charged particles that characterizes magnetic confinement, these models allow one to average over that fast gyration and reduce the 6D kinetic problem to a 5D gyrokinetic one.  That property is especially interesting for numerical simulation.  In this paper the finite Larmor radius approximation is derived from the Vlasov equation, in the limit of large uniform magnetic field and with an external electric field.  Because of the homogenization on the fast Larmor gyration, the limit equation \eqref{eq:Vlagyro} is  written in $5D$ gyro-coordinates $(x_g,v_{\parallel},|v_\perp|)$ defined in (\ref{eq:gyrocoor}). These coordinates are the position of the so-called particle guiding center, $x_g$ in the $3D$ space together with the parallel velocity $v_\parallel$ and the amplitude of transverse velocity $|v_\perp|$ that is proportional to the magnetic moment, an adiabatic invariant of the particle motion in the strong magnetic field limit (given a constant magnetic field) . Its mathematical structure is a combination of the Vlasov equation in the  direction parallel to the magnetic field and of the Euler equation in the perpendicular direction, where the original fields are replaced by the corresponding gyro-average fields. 

\medskip

To close the system, physicists use the electroneutrality equation $n_e=Z n_i$ where $n_e$ stands for the electron density and $n_i$ the ion density, each ion having a charge $Z$. In the following $Z=1$ will be considered with no loss of generality. It can be shown that the electroneutrality equation is in fact the the Poisson equation solved on scales that are significantly larger than the Debye length. Since the latter governs the Laplacian term of the Poisson equation, the electroneutrality equation is an appropriate approximation on scales of the order or larger than the Debye length.  Moreover, taking into account the difference between the density of particles and that of guiding centers leads one to introduce a polarization correction due to the non uniform distribution of particles on gyro-circles. In the second part of this paper, the electroneutrality equation \eqref{eq:elneut} is carefully written. The coupling of the 5D gyrokinetic Vlasov equation \eqref{eq:Vlagyro} to the electroneutrality equation  \eqref{eq:elneut}is the model used for instance in the GYSELA code, a project that aims at modeling the turbulent transport in fusion plasmas \cite{{Gysela06}, {Grandgirard_2006}}. While in the GYSELA code there is the possibility to use a  collision operator we concentrate here on the actual Vlasov equation with no collisions.
 
A difficulty raised by the electroneutrality equation taken as such is the lack of an explicit regularization term for the electric potential. Consequently, the regularity of the latter is not sufficient to ensure a mathematical solution of the Vlasov transport equation. 
The analysis of the well-posedness of limit model based on the electroneutrality equation thus seems impossible today. For that reason we restrict ourselves to the $1D$ Vlasov equation (without gyro-kinetic effects), coupled to the electroneutrality condition. That model is a particular case of the $3D$ model, provided the solutions do not depend on the direction perpendicular to the magnetic field.  Even for that simple kinetic model, not much is known about solutions. The difficulty lies in the fact that the force term is proportional to the derivative of the density in the field direction. Indeed, while the Vlasov equation ensures that all $L^p$-norms of the distribution function can be bounded, there is no control on the norms of its derivatives. In this paper, we prove the existence and uniqueness of a steady state solution in a slab geometry, therefore between two boundary conditions, provided one fulfills some conditions, in particular that they are no  particles trapped between the two boundaries.

\sectionnew{Derivation of the finite Larmor radius \\ approximation}

   The ion distribution in low density plasmas submitted to a  magnetic field is well described by the Vlasov equation. The latter is valid provided one can neglect the two-particle distribution function altogether \cite{Nicholson} so that the evolution of the standard distribution function is governed by the Liouville equation for the conservation of particles
\begin{equation} \label{eq:Vlasov}
\frac{\partial f}{\partial t} +   v \cdot \nabla_x f + \frac{Ze}{m_i}(E(t,x) + v \times B)\cdot \nabla_v f = 0~,
\end{equation}
where $f$ is the phase space density depending on time, position (in the domain $D$ of the plasma) and velocity (in $\R^3$) of the ions, and where $Ze$ and $m_i$ are  the ion charge and  mass respectively. As stated in the introduction, we shall consider $Z=1$ in the following with no loss of generality. A priori, the electric and magnetic fields are governed by the Maxwell equations, but a scale analysis allows one to approximate, and thus simplify, these laws. This will be made clear in the following.

   For a strong external magnetic field, the charged particles exhibit a fast rotation motion around the magnetic field lines. The frequency of that gyration, the Larmor frequency, is several orders of magnitude larger than the observed frequency range of turbulence. Furthermore, in the quasineutral limit this frequency is larger than all other frequencies of the particle dynamics, thus providing the means for an efficient scale separation. It is noteworthy that this is the basis of magnetic confinement that is presently realized in devices such as the tokamaks. In this framework, one separates two parts in the particle motion, on the one hand the slow motion of the center of the Larmor gyration, and, on the other hand, the fast gyro-motion. The phase space reduction achieved by only considering the slow motion, thus ignoring the fast gyro-motion, leads one to the $5D$ gyrokinetic model. These are usually derived by physicists either by averaging the single particle dynamics (\cite{Catto78} and \cite{Catto81}) or by a Lie-transform perturbative approach (\cite{Hahm88}, \cite{Brizard07} and \cite{Ghendrih08}).  Mathematically, this is achieved by looking at the limit of equation \eqref{eq:Vlasov} when the modulus of the external magnetic field $|B_{ext}|$ tends to infinity. However, different models can be obtained, depending on the way the other control parameters vary when $|B_{ext}|\rightarrow\infty$. For a magnetized plasma, two relevant scales characterize the limit regimes that have been discussed. These are the Debye length $\lambda_D$ and the thermal Larmor radius of the ions $\rho_{th}$. The Debye length weights the Laplacian in the Poisson equation and thus defines the transition from the non-neutral description of the plasma required on the sub-Debye scales from the quasineutral plasma description for scales larger than the Debye length. The thermal Larmor radius is the radius of the fast Larmor gyration for ions with average speed $v_{th}$. The averaging over the small time scales governed by the ion Larmor frequency then tends to translate into an averaging over the ion Larmor scale. The finite Larmor radius terms are therefore introduced to take into account this rather weak cut-off effect. The radius of the electrons gyration is generally much smaller for comparable ion and electron temperatures and will be neglected. The two scales discussed above are defined by
\begin{equation}
\lambda_D = \sqrt{\frac{\ep_0 T_e}{n e^2}}  \, , \qquad \rho_{th}= \frac{m_i v_{th}}{ e B} \, ,
\end{equation}
 
where the electron thermal energy $T_e$ is introduced rather than that of the ions since it is more appropriate to characterize the large electron mobility and therefore strong response to the electric field. 
Let us introduce the characteristic scale of the system $L$, for instance introduced by the boundary conditions. The large magnetic field limit then corresponds to the vanishing ion Larmor radius limit, namely the gyro-center approximation where $\rho_{th} / L \rightarrow 0$. Within this framework, E. Grenier uses pseudo-differential calculus to prove \cite{Gren97} the convergence of a 2D fluid model towards an incompressible fluid model when $\ep$ goes to zero with $\rho_{th} /L = \ep$ and $(\lambda_D / L)^2 = \ep$. For the same scaling and cold initial distributions, \emph{i.e.} the approximation of a Dirac measure at velocity zero for the distribution function, Y. Brenier proved in \cite{Bren00} with modulated energy technics the convergence of solutions of the Vlasov-Poisson system towards dissipative solutions (introduced by P.-L. Lions) of the Euler equation, for weel-prepared initial data. F. Golse and L. St-Raymond used the same technics in \cite{GolStR03} to prove the convergence of a Vlasov-Poisson model on the torus towards the $2D1/2$ Euler equation. This is done for a vanishing $\ep$ with $\rho_{th}/L = \ep$ and using another scaling for the Debye length, namely $\lambda_D / L= \ep$. This property is restricted to a torus size such that there is no resonant oscillation. This particular case appears to be relevant for electrons in a region close to the tokamak boundary.  

Derivations have also been performed  in the gyrocenter approximation, $\rho_{th} /L_\perp$ small but finite. Here $L_\perp$ is a carcteristic length in the perpendicular direction, which is usually choosen smaller thant the one in the parallel direction in finite Larmor radius approximation. In \cite{FreSon00} (see also \cite{FreSon-cras00} for a short version) E. Fr{\'e}nod and E. Sonnendr{\"u}cker studied the convergence of a linear Vlasov equation (external magnetic field), in this finite Larmor radius limit (large $B$ but finite Larmor radius), using two scale convergence methods. M. Bostan studied the 2D strong magnetic field limit, using Hilbert expansion technics, in a setting of finite Debye length and Larmor radius, $\lambda_D /(\ep L) = \rho_{th} / (\ep L) =1$, where $\ep$ is small but finite \cite{Bostan09} and obtained strong convergence results for regular initial conditions. Here, we further examine the linear case studied by Fr{\'e}nod and Sonnendr{\"u}cker, and obtain a simpler and more complete derivation. The finite Larmor radius approximation is correct for fusion plasmas, including the core and most of the edge plasma. However, the ion Larmor radius is much larger than the Debye length throughout the plasma. An interesting derivation should be done in the framework of present gyrokinetic calculations, namely $\rho_{th}/(\ep L) =1, \; \lambda_D /(\ep L)\rightarrow 0$, yet nothing has  been done in this very demanding limit (see the beginning of \ref{Sec:elecneut} for more details).

\medskip

The resolution of equation \eqref{eq:Vlasov} is known in the case where the fields $E$ and $B$ are external and $C^1$ or at least $BV$. In the $BV$-case, one can use the DiPerna Lions theory of transport equation \cite{DiPLio89} with its latest developments (\cite{Ambr04}). In the case where the electric potential is given by self-interaction with the usual Poisson law, and the magnetic field $B$ is still considered as external, results of existence and uniqueness obtained for the Vlasov equation without magnetic field may be used with the appropriate modifications. We refer to \cite{LioPer91}, \cite{Hors93}. In the case where the self-induced magnetic field is not negligible any more, we refer to the work \cite{DiPLio89-2} of DiPerna and Lions about the Vlasov-Maxwell equation. In the case of an electric field given by an electroneutrality equation, nothing is known from the mathematical point of view.

\subsection{The finite Larmor radius gyro-kinetic approximation. Rigorous approach.}
  
  For the sake of simplicity, we neglect the variation of the external magnetic field $B_{ext}$. In this cylindrical approximation, the curvature of the magnetic field is not considered, and one neglects the exploration of the magnetic field variation by the particles during the fast cyclotron motion. With respect to the relevant scales this effect if of the order of the Larmor radius $\rho_{th}$ divided by the characteristic radius of the magnetic field curvature, namely the tokamak major radius $R$. Although $\rho_{th} / R$ is a small parameter, this approximation is too strong with respect to many aspects of the tokamak physics, in particular since curvature effects are considered by physicists as the cause of a large class of micro-instabilities. The strong impact of such a small parameter can be traced back to a symmetry breaking of the system since the curvature governs the leading term that produces a charge separation. The present simplification must be regarded as a first step that is only valid when one addresses issues that do not lead one to symmetry breaking. Moreover, in order to avoid problems with boundary conditions, which are complex to handle, a periodic setting is used. In other words, we work on the torus $\T^3$ in space, and $\R^3$ in speed, with a constant magnetic field $B=|B|(0,0,1)$. 

  Next, we define dimensionless variables
\[
t'=\frac{t}{\tau}, \; x'_\para= \frac{x_\para}{L_\para}, \; x'_\perp = \frac{x_\perp}{L_\perp} , \; v'=\frac{v}{v_{th}} , \; E'=\frac{E}{E_0} \,,
\]
with characteristic time $\tau$, parallel length $L_\parallel$, perpendicular length $L_\perp$, velocity $v_{th}$ and electric field $E_0$. In this new set of dimensionless variables the Vlasov equation \eqref{eq:Vlasov} may be rewritten as:
\begin{equation} \label{eq:Vlares1}
\frac{\partial f}{\partial t'} +   \frac{\tau v_{th}}{L_\parallel}   \partial_{x'_\parallel} f + \frac{\tau v_{th}}{L_\perp} \cdot \nabla_{x'_\perp} f - \frac{\tau e E_0 }{m_i v_{th}}E \cdot \nabla_{v'} f - \frac{\tau e B}{m_i}v^\perp \cdot \nabla_{v'} f = 0 \,,
\end{equation}
where the subscript $\perp$ (resp. $\parallel$) stands for the projection on the perpendicular (resp. parallel) direction, and the superscript $\perp$ stands for the projection on the plane orthogonal to the field line and the rotation by $-\pi/2$ around the field line direction: $v^\perp=(v_2,-v_1,0)$ if $v=(v_1,v_2,v_3)$.  The time scale $\tau$ can conveniently be defined to reduce the number of control parameters in Eq.\eqref{eq:Vlares1} by setting:
\[
\frac{\tau v_{th}}{L_\parallel} =1 
\]
In a similar fashion, the normalization of the electric field can be used to define the third control parameter, hence:
\[
\frac{\tau e E_0}{m_i v_{th}} = 1
\]
The two remaining control parameters are then $L_{\parallel} / L_\perp$ for the third term in  Eq.\eqref{eq:Vlares1}   and $L_\parallel / \rho_{th}$ for the last term. It will be assumed here that both parameters $L_\perp$ and $\rho_{th}$ exhibit the same asymptotic behavior $L_\perp \propto \rho_{th}$. Let us define the perpendicular scale as $L_\perp= \rho_{th}  2 \pi /(k_\perp \rho_{th})$. In this expression, $k_\perp$ is the typical wave vector of the cross-field fluctuations. For a sufficiently small Larmor radius, one can then assume that the turbulence follows the so-called gyro-Bohm scaling such that $k_\perp \rho_{th}$ is a constant \cite{Lin}. In particular it does not depend on the magnitude of the magnetic field. One then finds the following relationship between the two control parameters
\[
\frac{L_\parallel}{\rho_{th}}= \frac{L_\parallel}{L_\perp} \frac{1}{k_\perp \rho_{th}}= \frac{1}{\ep} \,,
\] 
 that can be grasped by the following relation whenever one drops the proportionality constant
\[
\frac{L_\parallel}{\rho_{th}}= \frac{L_\parallel}{L_\perp} = \frac{1}{\ep} \,.
\] 
In practice, one could also define the normalization scale as $L_\perp= \rho_{th}$ so that the the two control parameters would be exactly the same. In this discussion, we have considered a transverse scale related to the fluctuation properties that allows one to characterize one of the terms in the Vlasov equation. However, other control parameters must be considered to account for the boundary conditions. The so-called $\rho_*$ parameter used in magnetic fusion is such a parameter since $\rho_*=\rho_{th}/a$ where $a$ is the plasma minor radius, the scale related to the boundary conditions while $\rho_{th}$ is the fluctuation scale. Relating the parameter $\ep$ to  $\rho_*$ introduces aspect ratios that stem from the periodic boundary conditions.
\[
\ep= \rho_*\frac{a}{L_\parallel}
\] 
The ratio $a / L_\parallel$ depends of the safety factor and the tokamak aspect ratio in the actual tokamak geometry but represents the ratio of the sizes of the domain in the radial and parallel direction with the present setting. The parameter $\ep$ can also be expressed in terms of the slow and fast characteristic times of the particle motion. The slow time  $\tau$ introduced in the normalization of Eq.\eqref{eq:Vlares1}, is the characteristic time to explore the tokamak geometry along the parallel direction and the fast time is due to the gyration motion, hence:
\[
\ep =\frac{1}{\tau\Omega}
\] 

In this last step one thus finds that the ion Larmor frequency $\Omega=e B/m_i$ is of order $(\tau\ep)^{-1}$. In other words, one thus assumes that the ion Larmor frequency is much larger than the parallel connection frequency $\tau^{-1}$. That is a valid assumption in all regions of a tokamak. Indeed, for ITER conditions, the ion Larmor frequency is of the order of $2~ 10^8 Hz$ and the connection time $\tau$ ranges from $10^{-3} s$ to $10^{-4} s$. These parameters mainly depend on the magnitude of the magnetic field and on the size of the device so that $\ep$ characterizes a given fusion device, $\ep \approx 10^{-5}$ for ITER parameters.

	Dropping the primes, we obtain the following rescaled version of the Vlasov equation \eqref{eq:Vlasov},
\begin{equation} \label{eq:Vlares}
\frac{\partial f}{\partial t} +   v_{_{\parallel}} \cdot \partial_{x_{_{\parallel}}} f + E \cdot \nabla_v f  +\frac{1}{\ep}(v_\perp \cdot \nabla_{x_\perp} f + v^\perp \cdot \nabla_{v_\perp} f )= 0~,
\end{equation}
 
\medskip  
  When $\ep \rightarrow 0$, the largest terms are those proportional to the factor $1/\ep$. Retaining the latter, we obtain a transport equation associated to the following system of ODE in the plane transverse to the magnetic field,
\begin{equation} \label{eq:ODE}
 \dot{x} = \frac{1}{\ep} ~v_{\perp } , \qquad \dot{v} = \frac{1}{\ep} v_{\perp }^\perp. 
\end{equation} 
As $B$ is homogeneous, the trajectories are circles of center $x_g = x_\perp + v_\perp^\perp$, and radius $|v_\perp|$ covered with frequency $\ep^{-1}$.
The global motion is the sum of this very quick motion of gyration and a slower and more complicated motion (with velocity of order one). In the limit of large $|B|$, particles are assumed to be evenly distributed on the gyro-circles, and their motion is the sum of a drift in the perpendicular direction, and a classical acceleration in the parallel direction. Heuristically, the electric drift $v_E$ may be obtained from the Newton law in normalized variables,
\begin{equation}
\dot{v} = E + \frac{1}{\ep}~ v^\perp  , \,
\end{equation}
Given the assumption that the particles have a fast motion of gyration, one can integrate the previous equation over a period of gyration, at lowest order, $v_E$ is given by the averaged velocity such that: 
\begin{equation}
0 \approx \langle E \rangle +   \frac{1}{\ep} ~v_E^\perp    \quad \text{so that }  v_E = \ep ~\langle E^\perp \rangle  \, .
\end{equation}
In the latter expression one readily recognizes the usual form of the electric drift velocity $v_E=E\times B / B^2$, however where  the electric field is averaged over a period $\langle E \rangle$. One can show that this average of the field $E$ corresponds to an average over a circle of radius $|v_\perp|$ in the perpendicular plane. Provided the only dependence on the gyrophase stems from the particle motion, this average translates into a Bessel operator defined as: 
\begin{equation} \label{eq:J0}
J^0_{\rho_{_L}} h(x_g) = \frac{1}{2\pi} \int_0^{2\pi} h(x_g + \rho_{_L} e^{i\varphi_c}) \,d\varphi_c \, ,
\end{equation}
 where $\rho_{_L}$ is the ion Larmor radius and \mbox{$e^{i\varphi_c} = (\cos\varphi_c,\sin\varphi_c,0)$}. We also introduce the operator $\tilde{J}^0_{\rho_{_L}}$ which stands for a position-velocity version of this average,
\begin{equation} \label{eq:J0tilde}
\tilde{J}^0_{\rho_{_L}} g(x_g,\rho_{_L},v_{\parallel}) = \frac{1}{2\pi} \int_0^{2\pi} g(x_g + \rho_{_L} e^{i\varphi_c}, \rho_{_L} e^{i(\varphi_c-\frac{\pi}{2})} + v_{\parallel} e_{\parallel}) \,d\varphi_c \, ,
\end{equation}
where $e_{\parallel} = (0,0,1)$. 

\medskip
Given the average on the gyrophase $\varphi_c$, the motion is reduced to a $5D$ space. The gyration radius $\rho_{_L}=|v_\perp|$ (given the chosen normalization) is related to the magnetic moment $\mu = m_i \rho_{_L}^2/ (2 e B)$, a quantity that is an adiabatic invariant, \emph{i.e.} a constant of motion in the large $B$ limit. One can show that this invariant is the conjugate variable of the gyrophase. Here, as $B$ is homogeneous, $\rho_{_L}=|v_\perp|$ will remain constant in the limit $|B_{ext}| \rightarrow +\infty$. We thus obtain a $4D+1D$ model, i.e. a $5D$ model with no dynamics in the variable derived from the magnetic moment. This reduction of the phase space is very interesting for numerical simulations. The following theorem states this reduction precisely.

\begin{theorem}
Let us assume that $f^0_\ep$ is uniformly bounded in $L^q$, $q>1$, that it weakly converges toward $f^0 \in L^q$, and that $E$ is a gradient and belongs to $L^p_{t,x}$ where $p^{-1}+q^{-1}=1$. For each $\ep>0$, let $f_\ep$ be a solution of \eqref{eq:Vlares} with the initial condition $f^0_\ep$. Then, up to a subsequence $(\ep_n)$, $(f_\ep(t,x_g-v^\perp,v))$ weakly converges to $\bar{f}$ in  $L^q$, where $\bar{f}$ only depends on \mbox{$(t,x_g,v_{\parallel},|v_\perp|=\rho_{_L})$} and is solution of
\begin{equation} \label{eq:Vlagyro}  
\frac{\partial \bar{f}}{\partial t} +   v_{\parallel} \, \partial_{x_{\parallel}} \bar{f} + J^0_{\rho_{_L}} E_{\parallel} \, \partial_{v_{\parallel}} \bar{f} + (J^0_{\rho_{_L}} E_\perp)^\perp \cdot \nabla_{x_g} \bar{f} = 0~,
\end{equation}
with the initial condition $\tilde{J}^0_{\rho_{_L}}(f^0)$.

Moreover, if $J^0_{\rho_{_L}} E \in BV(\R^3)$, for a.e. $\rho_{_L} >0$, then there is no need to extract a subsequence and the limit $\bar{f}$ is the unique solution of \eqref{eq:Vlagyro} with the initial condition $\tilde{J}^0_{\rho_{_L}}(f^0)$.

\end{theorem}

\begin{remark} 
The Bessel operator $J^0_{\rho_{_L}}$ has some regularization properties. It goes from $H^s$ to $H^{s+1/2}$ for all $s$, so that if $E \in H^1_{\parallel} \times H^{1/2}_\perp$, then $J^0_{\rho_{_L}} E \in H^1$, a condition that ensures the uniqueness of the solution of equation \eqref{eq:Vlagyro}.
\end{remark}

\begin{remark}
The operator $\tilde{J}^0$ is important to perform the adaptation of the 6D initial condition to the 5D limit model. Indeed, the very fast Larmor gyration creates an initial layer that instantaneously adapts the initial condition to the limit model.
\end{remark}

{\bf Proof of Theorem 1.1.}

The phase space in \emph{position-velocity coordinates} is not well adapted to perform the fast gyration averaging. To handle it more easily, it is convenient to change the system of coordinates and consider the \emph{gyro-coordinates} defined by
\begin{equation} \label{eq:gyrocoor}
 x_g = x + v^\perp, \qquad v_g=v \, .
\end{equation}
$x_g$ is the position of the gyro-center and $\rho_{_L}=|v^\perp|$ is the ion Larmor radius. To express the gradient in $x$ and $v$ in this new system of coordinates, one can conveniently remark that: 
\begin{equation}
\nabla_x h (dx_g - dv_g^\perp) + \nabla_v h dv_g = \nabla_{x_g} \bar{h} dx_g + \nabla_{v_g} \bar{h} dv_g,
\end{equation}
for any smooth function $h$, with the function $\bar{h}$ defined by $\bar{h}(x_g,v_g)=h(x,v)$. Then $\nabla_{x} = \nabla_{x_g}$ and \mbox{$\nabla_{v} = \nabla_{v_g} - \nabla_{x_g}^\perp$}, where \mbox{$\nabla_{x_g}^\perp = (\partial_{x_{g,2}},-\partial_{x_{g,1}},0)$}. Note that $\nabla^\perp$ stands for the gradient vector rotated by $-\pi/2$ and not $\pi/2$. The gradients are taken at the corresponding points. For instance, the first equality reads $\nabla_x h(x,v) = \nabla_{x_g} \bar{h}(x_g,v_g)$. The anti-symmetry of $\perp$, $a \cdot b^\perp = - a^\perp \cdot b$, has been used (and we shall make a wide use of it in the sequel). Given these relations,  \eqref{eq:Vlares} can be modified leading one to an equation satisfied by the function $\bar{f}_\ep(t,x_g,v_g) = f_\ep(t,x,v)$,
\begin{equation} \label{eq:Vlamod}
\begin{split}
\frac{\partial \bar{f}_\ep}{\partial t} +  & v_\parallel \, \partial_{x_\parallel} \bar{f}_\ep +   E_\parallel(t,x_g-v_g^\perp) \, \partial_{v_\parallel} \bar{f}_\ep \\
& + E_\perp(t,x_g-v_g^\perp) \cdot (\nabla_{v_g} \bar{f}_\ep - \nabla_{x_g}^\perp \bar{f}_\ep) +  \frac{1}{\ep}v_g^\perp \cdot \nabla_{v_g} \bar{f}_\ep= 0~,
\end{split}
\end{equation}
with initial condition $\bar{f}^0$. Here the subscript $\perp$ stands for the perpendicular components to the magnetic field, for instance $E_\perp = (E_1,E_2,0)$. Barred quantities are functions of the gyro-coordinates.

 Since (1.6) is a conservative transport equation, \emph{i.e.} it may be written as \mbox{$\partial_t f + \rm{div}(\dots) = 0$}, the $L^q$-norms of $f_\ep$ and then of $\bar{f}_\ep$ are conserved. Thus, $\| \bar{f_\ep} \|_{L^\infty_t(L^q_{x,v})}$ is uniformly bounded. Then, up to the extraction of a subsequence, we may assume that $\bar{f_\ep}$ weakly converges towards some $\bar{f} \in L^\infty_t(L^q_{x,v})$. Upon multiplying equation \eqref{eq:Vlamod} by $\ep$, in the limit $\ep\rightarrow 0$, we obtain
\[
v_{g}^{\perp} \cdot \nabla_{v_g} \bar{f} = 0 \, ,
\]
All the other terms from Eq.\eqref{eq:Vlamod} are bounded in the sense of distributions and their product with $\ep$ therefore vanishes in the limit $\ep\rightarrow 0$. This reduced form of the Vlasov equation implies that the only dependence of $\bar{f}$ on $v_{g,\perp}$ is on $|v_{g,\perp}|$. Hence with no dependence on the gyrophase.

Let us now consider equation \eqref{eq:Vlamod} for $\bar{f}_\ep$, when integrated against a smooth test fonction $\phi$ with support in time avoiding $t=0$ (we will handle the initial conditions later). Let us further assume that the dependence on $v_{g,\perp}$ is restricted to a dependence on $|v_{g,\perp}|$. For such a function, one readily finds that
\[
\frac{1}{\ep}\int \bar{f}_\ep v_g^\perp \nabla \phi dx_gdv_g =0
\]
for symmetry reasons. As a consequence, the projection of $\bar{f}_\ep$ takes the following form
\begin{equation} \label{eq:vladist}
\int \bar{f}_\ep (\partial_t \phi + v_{g,\parallel} \partial_{x_{g,\parallel}} \phi + E_\parallel \partial_{v_{g,\parallel|}} \phi + E_\perp \cdot \nabla_{v_g} \phi - E^\perp \cdot \nabla_{x_g} \phi ) \,dx_g dv_g dt = 0~,
\end{equation}
keeping in mind that $E$ is calculated at the point $x_g-v_g^\perp$. We have also used the relation $a^\perp \cdot b = - a \cdot b^\perp$ as well as the fact that equation \eqref{eq:Vlamod} may be written in a conservative form because ( if $J$ is the matrix of the linear map $v \rightarrow v^\perp$),
\begin{eqnarray*}
\diver_{v_g} (E(t,x-v_g^\perp)) & = & \Tr ((\nabla E) \, J) = - \partial_{x_2} E_1 + \partial_{x_1} E_2 = 0~, \\
\diver_{x_g} E^\perp & = & \Tr (J \, (\nabla E)) = \Tr((\nabla E) \, J) = 0 \, ,
\end{eqnarray*}
since $E$ is a gradient. At this stage, we can take the limit and obtain that $\bar{f}$ also satisfies equation \eqref{eq:vladist}.

This is not yet a proper equation in the sense of distributions, even though $f$ depends on $|v_{g,\perp}|$. Indeed,  the electric field still exhibits the dependence on the particle position $E\equiv E(t,x_g - v_g^\perp)$. To obtain an equation only depending on $|v_{g,\perp}|$, we use polar coordinates for $v_{g,\perp }$, \emph{i.e.} $v_g = ( \rho_{_L} \cos \varphi_c, \rho_{_L} \sin \varphi_c,v_\parallel)$ and Fubini's theorem to integrate the previous integral first in $\varphi_c$, then in the other variables. This leads one to the following equation:
\begin{equation} \label{findemo}
\begin{split}
\int_{x_g,v_{g,\parallel},r} (2\pi \rho_{_L}) \bar{f} (\partial_t \phi + & v_{g,\parallel} \partial_{x_{g,\parallel}} \phi + J^0_{\rho_{_L}} E_\parallel \partial_{v_{g,\parallel}} \phi + (J^0_{\rho_{_L}} E)^\perp \cdot \nabla_{x_g} \phi \\
&  + R(t,x_g,v_{g,\parallel},\rho_{_L})\cdot\nabla _{\rho_{_L}}\phi )\,dx_g dv_{g,\parallel}  d\rho_{_L} =0~,
\end{split} 
\end{equation}
where the the operator $J^0_{\rho_{_L}}$ is defined in \eqref{eq:J0}, and the term $R$ is equal to
\[
R(t,x_g,v_{g,\parallel},\rho_{_L}) =  \int_0^{2\pi} E_\perp(x_g - \rho_{_L} e^{i(\varphi_c+ \frac{\pi}{2})}) \cdot e^{i\varphi_c} \,d\varphi_c   \, ,
\]
with the previously introduced notation $e^{i\varphi_c} = (\cos\varphi_c,\sin\varphi_c,0)$. It can be shown that $R$ is null since it is the circulation of the electric field $E$ along a gyrocircle. The latter vanishes in the electrostatic case since $E$ is a gradient. With $R=0$, \eqref{findemo} is identical to the equation \eqref{eq:Vlagyro} written distribution wise for $2\pi \rho_{_L} \bar{f}$. As there is no dynamics in the $\rho_{_L}$ direction, the factor $\rho_{_L}$ is only a multiplicative constant. One can then introduce $\scriptstyle \Phi /\sqrt{\rho_{_L}^2 + \eta}$ as test function, and then obtain that $\scriptstyle \bar{f} \rho_{_L}/\sqrt{\rho_{_L}^2 + \eta}$ satisfies the equation. Letting $\eta \rightarrow 0$ in that linear equation, we obtain the appropriate result for the equation governing the evolution of $\bar{f}$.

Regarding the initial conditions we use a similar projection technique with a test function $\Phi$ depending on $|v_{g,_{\parallel}}|$, but  not vanishing at $t=0$. We then obtain equation \eqref{eq:vladist}, with the right-hand side replaced by
\begin{equation}
- \int \bar{f_\ep}^0 \Phi(0,x_g,v_{g,\parallel},\rho_{_L}) \, dx_g dv_{g,\parallel}  \rho_{_L} d\rho_{_L} d\varphi_c \, .
\end{equation}
The asymptotic limit then leads one to:
\begin{equation}
\int \frac{1}{2\pi} \left(\int_0^{2\pi} \bar{f}^0(x_g,\rho_{_L} e^{i\theta}+ v_\parallel e_{\parallel})  d\varphi_c \right) \Phi(0,\cdot) \, dx_g dv_{g,\parallel} 2\pi \rho_{_L} \, d\rho_{_L} \, .
\end{equation}
Changing coordinate to recover the \emph{position-velocity} coordinates in the $\varphi_c$ integral, we exactly obtain the initial condition $\tilde{J}_{\rho_{_L}} f^0$ that is expected for $\bar{f}$. 
\cqfd

This proof is valid for external electric and magnetic fields. The case where the equation \eqref{eq:Vlamod}, invariant in the parallel direction, is coupled to the Poisson equation given a Debye length of the same order as the Larmor radius 
has been treated by Fr{\'e}nod and Sonnendr{\"u}cker in \cite{FreSon01}. It corresponds to a Debye length of order $\sqrt{\ep}$. Technically, we may handle that case with our technic to obtain weak solution as in Arsenev's work \cite{Arsenev75}. The main point is to obtain $L^p$ estimates on the density $n_i = \int f \,dv$. They can be obtained by classical estimates using upper bounds for the kinetic energy. The case where the Debye length is taken much smaller than the Larmor radius will be of greater interest, but the difficult problem is there to average the strong oscillations appearing at the scale of the Debye length. We refer to \cite{Gren95}, \cite{Gren96} and \cite{CorGre00} for more details on quasi-neutral plasma without magnetic field, and to \cite{GolStR03} for results in the guiding-center approximation with $\rho_{_L}\propto\ep$, $\lambda_D\propto \ep$ and $\ep \rightarrow 0$.

\subsection{The electroneutrality equation. Heuristic approach.} \label{Sec:elecneut}

In this section we address the self-consistent problem when linking the electric field to the charge distribution. The relevant equation is the Maxwell-Gauss equation that relates the divergence of the electric field to to the local charge governed by the particle density of charged particles. The latter must be determined using the ion distribution function for the guiding centers that is solution of the gyrokinetic Vlasov equation. A similar treatment for the electrons must be done. We will follow  heuristic arguments together with assumptions that are not justified rigorously. However, this heuristic derivation bares some interest. It is an alternative to the physicists' presentation of that equation based on the Fourier transform. Furthermore, it allows one to recover the electroneutrality equation \eqref{eq:elneut} used in the GYSELA code to close the gyrokinetic equation \eqref{eq:Vlagyro}.

\medskip

For quasineutral plasmas, the electric field response to any charge separation
 governs a restoring force. Should the charge separation extend on a scale larger than the Debye length, the restoring force would be to strong to allow any significant charge build-up. The plasma can thus be considered to be everywhere with near zero charge, hence quasineutral. As a consequence, the
density of negative charges $e n_e$ ($n_e$ being the density of
electrons) equals the density of positive charges $e n_i$ ($n_i$ being
the density of ions). In the electrostatic limit, neglecting the time dependence on the vector potential, one can recover this physics based argument as an asymptotic limit of the Poisson equation for the 
electric potential. 
\[ 
\frac{e}{T_e} \lambda_D^2 \Delta \Phi = \frac{n_e- n_i}{n_e} \, ,
\] 
The right hand side is dimensionless and so is $e\Phi/T_e$. The dimensionless control parameter on the left hand side operator thus appears as the square of the ratio of the Debye scale divided by the characteristic scale of the charge separation. In the limit $\lambda_D^2 \rightarrow 0$,  one readily recovers the quasineutrality equation, namely: $n_e = n_i$. Let us consider the ions density. Given the Vlasov equation and its dependence on the electric field, it is likely that the ion distribution function (and thus the ion density) is an implicit function of the electrostatic potential. The same applies for the electron density. Then, in the limit $\lambda_D^2 \rightarrow 0$,  there is no analytical dependence on the electric potential. The latter then becomes a Lagrangian multiplier associated to the quasineutrality equation.   
On such issues, we acknowledge the work of E. Grenier \cite{Gren95},
\cite{Gren96} (the only work, at our knowledge, on that subject where the limiting
model is kinetic and not fluid), Y. Brenier \cite{BreGre94}, S. Cordier
\cite{CorGre00} and N. Masmoudi \cite{Masm01}. 
 
However, if we assume that the system is close to equilibrium, \emph{i.e.} that the departure from a constant electric potential is small, then 
the electroneutrality equation $n_e = n_i$ provides an explicit dependence 
on a mean-field electric potential. For an electron-ion plasma, the mass ratio is such that the electrons are far more mobile than the ions. One can then
assume that their response to an electrostatic perturbation
is adiabatic on magnetic field lines or surfaces. In that
case, assuming that the equilibrium density of electrons $n_{e,0}$ and that of 
ions $n_{i,0}$ are equal to $n_0$, one may write
\begin{equation}
n_e = n_0 e^{^{\frac{e}{T_e}(\Phi - \langle \Phi \rangle )}} \approx n_0 \big(1 + \frac{e}{T_e}( \Phi - \langle \Phi \rangle ) \big) \, ,
\end{equation}
where $\langle \Phi \rangle$ is the average of $\Phi$ on a closed
magnetic field line or surface.  The expansion on the right hand side holds if $e\Phi << T_e$. In the very simple geometry that is considered here (with $\mathbf{B}=(0,0,B)$), $\langle \Phi \rangle $ is the average in the $x_{\parallel}$ direction: $\langle \Phi \rangle = \int \Phi(x) \,dx_{\parallel}$. The approximation of adiabatic electrons thus reintroduces an explicit dependence on the electrostatic potential.

For the perturbation of the ion distribution, the first difficulty is
to obtain the distribution of ions in physical space from
$\bar{f}$ written in gyro-coordinates. If the distribution in physical
space is assumed to be constant on gyrocircles as shown in the last
section for the limit model, the following formula is obtained,
\begin{equation}
f(t,x,v) =\bar{f}\big(t, x +  v^\perp, |v_\perp|, v_{\parallel} \big) \, . 
\end{equation}
Taking the integral in $v$ leads to 
\begin{equation}
n_i(t,x) = n_0\int J_{\rho_{_L}}^0\big(\bar{f}(t, x,\rho_{_L},v_{\parallel})\big) 2 \pi \rho_{_L} \,d\rho_{_L} dv_{\parallel} \, ,
\end{equation}
where $J_{\rho_{_L}}^0$ only acts on the $x$ variable (since $\rho_{_L}=|v_\perp|$ with our conventions). The quasineutrality equation, $n_e=n_i$, then leads one to
\begin{equation}
1 + \frac{e}{T_e} \big( \Phi - \langle \Phi \rangle \big) =  \int
J_{\rho_{_L}}^0\big(\bar{f}(t, x,\rho_{_L},v_{\parallel})\big) 2 \pi \rho_{_L} d\rho_{_L} dv_{\parallel} .
\end{equation} 
However, the assumption that the distribution $f$ is constant
on the gyrocircles is unrealistic whenever the electric potential is not constant. In order to express the inhomogeneity of the density on gyrocircles, we  add a perturbation to $\bar{f}$. It can be expressed as an adiabatic perturbation on the gyrocircles
\begin{equation}
- \frac{e}{T_i} (\Phi - \bar{\Phi}) n_0 f_i(v) \, ,
\end{equation}
$\bar{\Phi}$ is the average of $\Phi$ over the gyrocircle of a given particle, and $f_i(v)$ is the  equilibrium distribution of ions. 
Note that this adiabatic perturbation depends on $v_\perp$
through the choice of the gyrocircle used in determining the average $\bar{\Phi}$.
\begin{eqnarray}
{\textstyle \bar{\Phi}(t,x,v_\perp) } & = & \frac{1}{2\pi} {\textstyle \int \Phi \big( t, x +  v^\perp + \rho_{_L} e^{i\varphi_c} \big) \,d\varphi_c \, , } \nonumber \\
n_0 {\textstyle \int \bar{\Phi}(t,x,v_\perp) f_i(v)\,dv } & = & \frac{n_0}{2\pi}
{\textstyle \int \Phi \big(t,x -  \rho_{_L} (e^{i\varphi_c'} + e^{i\varphi_c}) \big) f_i(\rho_{_L},v_{\parallel}) \,d\varphi_c d\varphi_c' \rho_{_L} d\rho_{_L} dv_{\parallel} } \nonumber\\
& = &  n_0{\int (J_{\rho_{_L}}^0)^2 \Phi (t,x) h_i(\rho_{_L})\,d\rho_{_L} \, ,}
\end{eqnarray}
where $e^{i \varphi_c} =(\cos \varphi_c,\sin \varphi_c,0) $ and $h_i(\rho_{_L})
= 2 \pi \rho_{_L} \int f_i(\rho_{_L},v_{\parallel}) \,dv_{\parallel}$. When the equilibrium distribution is a maxwellian,
$f_i(v) = 1 /(T_i \pi )^{3/2} e^{-|v|^2/T_i}$, so that $h_i(\rho_{_L})
= 2 \rho_{_L}/T_i  \,e^{-  \rho_{_L}^2/T_i}$.
Adding the adiabatic perturbation of the gyrocircles, we finally obtain the following electroneutrality equation,
\begin{equation} \begin{split}
1 + \frac{e}{T_e}\big(\Phi - \langle \Phi \rangle\big)  =   \int
J^0_{\rho_{_L}} & \big(\bar{f}(t, x,\rho_{_L},v_\parallel)\big) 2 \pi\,\rho_{_L}
d\rho_{_L} dv_\parallel  \\ 
- & \frac{e}{T_i}  \int \big(1 -(J^0_{\rho_{_L}})^2\big) \Phi (t,x)
h_i(\rho_{_L}) d\rho_{_L} \, .
\end{split} \end{equation}
 Multiplying by $T_e/e$, then yields
\begin{equation} \label{eq:elneut}
\begin{split}
\big(\Phi - \langle \Phi \rangle \big) +  & \frac{T_e}{T_i} \int \big( \Phi - 
(J^0_{\rho_{_L}})^2 \Phi \big)\, h_i(\rho_{_L}) d\rho_{_L}  =  \\
& \frac{ T_e}{e} \Big( \int J^0_{\rho_{_L}}\big(\bar{f}(t, x,\rho_{_L},v_\parallel)\big) 2
  \pi \rho_{_L}d\rho_{_L} dv_\parallel - 1\Big) \, .
\end{split}
\end{equation}

\begin{remark} \label{rem:convol}
If the equilibrium density of ions is a Maxwellian, then  $h_i(\rho_{_L})= 2 \rho_{_L}/T_i  \,e^{-
  \rho_{_L}^2/T_i}$ and the operator 
\begin{equation}
\int (J^0_{\rho_{_L}})^2  \, h_i(\rho_{_L}) d\rho_{_L}
\end{equation}
is the convolution in the perpendicular plane with the
radial function $H_{_{T_i}}(r)$ defined by
\begin{equation}
H_{_{T_i}}(r) = \frac{e^{-\frac{r^2}{4 T_i}}}{2 \pi^{3/2} \sqrt{T_i} r}.
\end{equation}
See Appendix \ref{App:A} for a detailed proof.
\end{remark}

Since $\rho_{_L}$ is a parameter in the equation of motion \eqref{eq:Vlagyro}, the equilibrium distribution and the initial perturbation remain concentrated on a unique value of $\rho_{_L}$ at all time provided it is the case for the initial conditions. A first step to solve the system might be to start under this assumption of single value of  $\rho_{_L}$. This would likely be the most difficult step, since the general case is a superposition of such cases. The gyrokinetic Vlasov equation \eqref{eq:Vlagyro} is unchanged when considering a single value of $\rho_{_L}$, however the electroneutrality equation \eqref{eq:elneut} may be simplified and becomes:
\begin{equation} \label{eq:elneutsim}
\big(\Phi - \langle \Phi \big) +  \frac{ T_e}{T_i} \Big(1 -
(J^0_{\rho_{_L}})^2\Big) \Phi (t,x)  = \frac{ T_e}{e} \Big(  J^0_{\rho_{_L}}\big(\bar{n_i}(t, x)\big) - 1 \Big) ~,
\end{equation} 
where $\bar{n_i}(t,x) = 2 \pi \rho_{_L} \int \bar{f}(t,x,v_\parallel) \,dv_\parallel$ is the density of ions in gyro-coordinates.
 
Still, after this further simplification, the gyrokinetic Vlasov equation \eqref{eq:Vlagyro} coupled to the electroneutrality equation \eqref{eq:elneutsim} remains is a very difficult mathematical problem. The main difficulties are the following.
\begin{itemize}
\item The lack of regularity in the parallel direction. The potential $\Phi$ has the regularity of $f$ in the parallel direction. This is to be compared to the Vlasov-Poisson case where $D^2 \Phi$ ($\Delta \Phi$ in the Poisson equation given above) has the regularity of $f$. We may overcome this problem by adding some viscosity in that direction, in other words by adding a term $- \lambda \partial_{v_\parallel}^2 \bar{f}$ in equation \eqref{eq:Vlagyro}. 
\item A less important lack of regularity lies in the perpendicular direction. In fact, only the gyro-average of $\Phi$ appears in \eqref{eq:Vlagyro}. Moreover, the term $\langle \Phi \rangle $ is more regular than $\Phi $, so that by 
\eqref{eq:elneutsim} $\Phi$ has the regularity of $J^0_{\rho_{_L}}(\bar{\rho})$. Hence $J^0_{\rho_{_L}}(\Phi)$ has the regularity of $(J^0_{\rho_{_L}})^2(\bar{\rho})$. The operator $(J^0_{\rho_{_L}})^2$ sends $L^2$ into $H^1$. That is ``almost'' enough, in the sense that we could reach a sufficient regularity provided $(J^0_{\rho_{_L}})^2$ were compact from $L^2$ in $H^1$.
\end{itemize}

In view of these difficulties, we focus in the next section on the lack of regularity in the parallel direction, and consider time-independent solutions depending only on $x_\parallel$ and $v$, of the form $f(x_\parallel,v_\parallel) f_\perp(|v_\perp|)$.

\sectionnew{Steady state solutions in the direction parallel to the magnetic field}
Let us consider steady state solutions $(f,\Phi )$ to \eqref{eq:Vlagyro}-\eqref{eq:elneut} and let us assume that the function $f$ can be written in the form $f(x_{\parallel},v_{\parallel}) f_\perp(|v_\perp|)$, with $\int_0^{+\infty} f_\perp(|v_\perp|) 2 \pi |v_{\perp}| \,d|v_{\perp}| =1$. Then the  term $f_i$ has no incidence in the evolution equation and can thus be factorized in \eqref{eq:Vlagyro}. Taking $T_e= 1$ for the sake of simplicity, taking again into account the equilibrium density $n_0$, not always constant in that section and therefore denoted by $n$, and replacing $n_i$ by $\rho$, the electroneutrality equation \eqref{eq:elneut} then reads: 
\[
\Phi -<\Phi >= \frac{\rho }{n}-1,
\]
where $\rho (z)= \int f(z,v)dv$. To further simplify the notations, let $z=x_{\parallel}$ and $v=v_{\parallel}$. Solving \eqref{eq:Vlagyro}-\eqref{eq:elneut} is then equivalent to finding a distribution function $f$ solution of the following equation. 
\begin{equation} \label{eq:Vlas_red}
v\frac{\partial f}{\partial z}-\Big( \frac{\rho }{n}\Big) '\frac{\partial f}{\partial v}= 0 \, 
\end{equation}
Note that $g '$ is the derivative of $g$ with respect to $z$. Furthermore, we consider the problem in a slab geometry $z\in [-1,1]$, and therefore require given $f_{\pm }$ as boundary conditions,
\begin{equation}\label{bc:Vlas_red}
f(-1,v)= f_-(v),\hspace*{0.05in}v>0,\quad f(1,v)= f_+(v),\hspace*{0.05in}v<0.
\end{equation} 
\begin{lemma}
Given $n$ positive, then any solution $f$ to \eqref{eq:Vlas_red}-\eqref{bc:Vlas_red}, such that $\frac{\rho }{n}$ is non decreasing, must satisfy:
\begin{eqnarray} \label{eq:solution_rho}
{ \textstyle \rho (z)= \int _{-\infty }^0\frac{\mid u\mid f_+(u)}{\sqrt{u^2+2\frac{\rho }{n}(1)- 2\frac{\rho }{n}(z)}}du+2\int _{\scriptscriptstyle{\sqrt{2\frac{\rho }{n}(z)-2\frac{\rho }{n}(-1)}}}^{\scriptscriptstyle{\sqrt{2\frac{\rho }{n}(1)-2\frac{\rho }{n}(-1)}}}\frac{uf_-(u)}{\sqrt{u^2-2\frac{\rho }{n}(z)+2\frac{\rho }{n}(-1)}}du\nonumber } \\
{ \textstyle + \int _{\scriptscriptstyle{\sqrt{2\frac{\rho }{n}(1)-2\frac{\rho }{n}(-1)}}}^{+\infty }\frac{uf_-(u)}{\sqrt{u^2-2\frac{\rho }{n}(z)+2\frac{\rho }{n}(-1)}}\,du \, ~~~~~~ }
\end{eqnarray}
\end{lemma}
\underline{Proof of Lemma 2.1}\hspace*{0.05in}The characteristic $(Z,V)$ of \eqref{eq:Vlas_red} starting from $(z,v)$ is defined by:
\begin{eqnarray*}
Z'(s)= V(s),\quad Z(0)= z,\\
V'(s)= -\Big( \frac{\rho }{n}\Big) '(Z(s)),\quad V(0)= v.
\end{eqnarray*}
Hence
\[
V^2(s)+2\frac{\rho }{n}(Z(s))= v^2+2\frac{\rho }{n}(z).
\]
For $v>0$, it crosses $\{ (-1,u), u>0\} $. For $v<0$, it crosses $\{ (1,u), u<0\} $ if and only if there is a solution $V$ to $V^2+2\frac{\rho }{n}(1)= v^2+2\frac{\rho }{n}(z)$, i.e. $v^2\geq 2\Big( \frac{\rho }{n}(1)-\frac{\rho }{n}(z)\Big) $. Consequently,
\begin{eqnarray*}
f(z,v)= f_-(V(s_-(z,v)))\hspace*{0.08in} if \hspace*{0.08in}v>-\sqrt{2\Big( \frac{\rho }{n}(1)-\frac{\rho }{n}(z)\Big) },\\
f(z,v)= f_+(V(s_+(z,v)))\hspace*{0.08in} if \hspace*{0.08in}v<-\sqrt{2\Big( \frac{\rho }{n}(1)-\frac{\rho }{n}(z)\Big) },
\end{eqnarray*}
where 
\[
V^2(s_{\pm }(z,v))+2\frac{\rho }{n}(\pm 1)= v^2+2\frac{\rho }{n}(z),\quad V(s_-(z,v))>0, \quad V(s_+(z,v))<0.
\]
Hence
\begin{eqnarray*}
f(z,v)= f_-\Big( \sqrt{2\Big( \frac{\rho }{n}(z)-\frac{\rho }{n}(-1)\Big) +v^2}\Big) \hspace*{0.07in} if \hspace*{0.07in}v>-\sqrt{2\Big( \frac{\rho }{n}(1)-\frac{\rho }{n}(z)\Big) },\\
f(z,v)= f_+\Big( -\sqrt{2\Big( \frac{\rho }{n}(z)-\frac{\rho }{n}(1)\Big) +v^2}\Big) \hspace*{0.08in} if \hspace*{0.08in}v<-\sqrt{2\Big( \frac{\rho }{n}(1)-\frac{\rho }{n}(z)\Big) }.
\end{eqnarray*}
Consequently,
\begin{eqnarray*}
\rho (z)= \int _{-\infty }^{-\sqrt{2(\frac{\rho }{n}(1)-\frac{\rho }{n}(z))}}f_+\Big( -\sqrt{2(\frac{\rho }{n}(z)-\frac{\rho }{n}(1))+v^2}\Big) dv\\
+\int _{-\sqrt{2(\frac{\rho }{n}(1)-\frac{\rho }{n}(z))}}^{+\infty }f_-\Big( \sqrt{2(\frac{\rho }{n}(z)-\frac{\rho }{n}(-1))+v^2}\Big) dv.
\end{eqnarray*}
Changes of variables in both integrals lead to \eqref{eq:solution_rho}.\\
For a constant density $n$, trivial solutions to \eqref{eq:Vlas_red}-\eqref{bc:Vlas_red} are 
\[ f(z,v)= f_-(v),\hspace*{0.05in}v>0,\quad f(z,v)= f_+(v),\hspace*{0.05in}v<0. \]
Proving the existence of solutions to \eqref{eq:Vlas_red}-\eqref{bc:Vlas_red} satisfying $\Big( \frac{\rho }{n}\Big) '\geq 0$, for a non constant $n$, is the aim of this section. Denote by $H_i$, $1\leq i\leq 4$, the following set of assumptions.
\begin{eqnarray*}
(H_1)\quad n'\leq 0.\\
(H_2)\quad \sup_{x\geq 0}\int _{\sqrt{x}}^{+\infty }\frac{uf_-(u)}{\sqrt{u^2-x}}du=: \lambda <+\infty .\\
(H_3)\quad f_-(u)= 0,\quad 0<u<2\sqrt{\frac{2\mu}{n(1)}},\hspace*{0.05in} \text{where } \hspace*{0.05in}\mu = \int _{-\infty }^0f_+(u)du+\lambda .\\
(H_4)\quad \int _0^{+\infty }f_-(u)du<\int _{-\infty }^0f_+(u)du,\quad \int _{-\infty }^0\frac{f_+(u)}{u^2}du<\frac{n(1)}{4}.
\end{eqnarray*}
Notice that $H_2$ is satisfied when 
\begin{eqnarray*}
\lim _{u\rightarrow +\infty }uf_-(u)<+\infty \hspace*{0.07in} and \hspace*{0.05in} \int _0^{+\infty }u\mid f_-'(u)\mid du<+\infty .
\end{eqnarray*}

\begin{remark}
The assumption $(H2)$ ensures that $\rho$ is a perturbation of the equilibrium $n$, so that with assumption $(H1)$ the force-field will always be oriented rightward (no possibilities of trapped particles). \\
The assumption $(H3)$ ensures that there are two beams of ions (one coming from the right and the other from the left), and that all the ions coming from one side reach the other side (no turn-back).\\
The $(H4)$ assumption is more technical and enforces that the derivative of the density is bounded. 
\end{remark}

\begin{lemma}
Assume $H_i$, $1\leq i\leq 4$, and $\frac{n'}{n^2}\in L^{\infty }$. Denote by
\begin{eqnarray}
{ \textstyle
K= \{ \alpha \in W^{1,\infty }([-1,1]); \alpha \geq 0, \alpha (1)-\alpha (-1)\leq \frac{2\mu }{n(1)}, 0\leq \alpha '\leq 4\mu \parallel \frac{n'}{n^2}\parallel _{\infty }\} \, . ~~}
\end{eqnarray}
There is a solution $\alpha \in K$ to 
\begin{eqnarray} \label{eq:alpha}
\frac{n}{2}\alpha (z)= \int _{-\infty }^0\frac{\mid u\mid f_+(u)}{\sqrt{u^2+\alpha (1)-\alpha (z)}}du+\int _{2\sqrt{\frac{2\mu }{n(1)}}}^{+\infty }\frac{uf_-(u)}{\sqrt{u^2-\alpha (z)+\alpha (-1)}}du  \, ,\nonumber \\
z\in [-1,1]\hspace*{0.05in}.~~
\end{eqnarray}
Moreover, $\alpha $ is the unique non-decreasing solution of \eqref{eq:alpha} such that 
\begin{eqnarray*}
\alpha (1)-\alpha (-1)\in [0, \frac{2\mu }{n(1)}].
\end{eqnarray*}
\end{lemma}
\underline{Proof of Lemma 2.2.}\hspace*{0.05in}Prove that the map $F$ that maps $\alpha \in K$ in $\beta $ defined by
\begin{eqnarray*}
\frac{n}{2}\beta (z)= \int _{-\infty }^0\frac{\mid u\mid f_+(u)}{\sqrt{u^2+\alpha (1)-\alpha (z)}}du\\
+\int _{2\sqrt{\frac{2\mu }{n(1)}}}^{+\infty }\frac{uf_-(u)}{\sqrt{u^2-\alpha (z)+\alpha (-1)}}du, \quad z\in [-1,1],
\end{eqnarray*}
has a fixed point. First, $F$ maps $K$ in $K$. Indeed, $\beta \in W^{1,\infty }([-1,1])$ like $\alpha $, \\
$\beta $ is nonnegative, and
\begin{eqnarray*}
\frac{n(1)}{2}\beta (1)= \int _{-\infty }^0f_+(u)du+\int _{\sqrt{\alpha (1)-\alpha (-1)}}^{+\infty }\frac{uf_-(u)}{\sqrt{u^2-\alpha (z)+\alpha (-1)}}du\\
\leq \int _{-\infty }^0f_+(u)du+\int _{\sqrt{\alpha (1)-\alpha (-1)}}^{+\infty }\frac{uf_-(u)}{\sqrt{u^2-\alpha (1)+\alpha (-1)}}du\leq \mu .
\end{eqnarray*}
Hence, $\beta (1)-\beta (-1)\leq \frac{2\mu }{n(1)}$. Moreover,
\[
\beta '= 2\frac{\mid n'\mid }{n^2}X+\frac{\alpha '}{n}Y,
\]
where
\begin{eqnarray*}
X= \int _{-\infty }^0\frac{\mid u\mid f_+(u)}{\sqrt{u^2+\alpha (1)-\alpha (z)}}du+\int _{2\sqrt{\frac{2\mu }{n(1)}}}^{+\infty }\frac{uf_-(u)}{\sqrt{u^2-\alpha (z)+\alpha (-1)}}du,\\
Y= \int _{-\infty }^0\frac{\mid u\mid f_+(u)}{(u^2+\alpha (1)-\alpha (z))^{\frac{3}{2}}}du+\int _{2\sqrt{\frac{2\mu }{n(1)}}}^{+\infty }\frac{uf_-(u)}{(u^2-\alpha (z)+\alpha (-1))^{\frac{3}{2}}}du.
\end{eqnarray*}
Since $\beta '\geq 0$, and, $X\leq \int _{-\infty }^0f_+(u)du+\lambda = \mu $, then,
\begin{eqnarray*}
\sup_{x\geq 0}\int _{2\sqrt{x}}^{+\infty }\frac{u}{(u^2-x)^{\frac{3}{2}}}f_-(u)du\leq \frac{n(1)}{3\sqrt{3}\mu }\int _0^{+\infty }f_-(u)du.
\end{eqnarray*}
Indeed, either $x\geq \frac{2\mu }{n(1)}$ and then
\begin{eqnarray*}
u>2\sqrt{x} \Rightarrow \frac{u}{(u^2-x)^{\frac{3}{2}}}\leq \frac{2}{3\sqrt{3}x}\leq \frac{n(1)}{3\sqrt{3}\mu },
\end{eqnarray*}
or $x<\frac{2\mu }{n(1)}$ and then $f_-(u)\neq 0$ implies that $u>2\sqrt{\frac{2\mu }{n(1)}}$, hence
\begin{eqnarray*}
\frac{u}{(u^2-x)^{\frac{3}{2}}}\leq 
\frac{2\sqrt{\frac{2\mu }{n(1)}}}{(\frac{8\mu }{n(1)}-x)^{\frac{3}{2}}}
\leq \frac{n(1)}{3\sqrt{3}\mu }.
\end{eqnarray*}
So that, 
\begin{eqnarray*} \label{eq:pourG}
\sup _{x\geq 0}\int _{2\sqrt{x}}^{+\infty }\frac{u}{(u^2-x)^{\frac{3}{2}}}f_-(u)du\leq \frac{n(1)}{3\sqrt{3}\int _{-\infty }^0f_+(u)du}<\frac{n(1)}{4},
\end{eqnarray*}
and by $H_4$,
\[
Y\leq \int _{-\infty }^0\frac{f_+(u)}{u^2}du+\frac{n(1)}{4}\leq \frac{n(1)}{2}.
\]
Consequently, $\beta '\leq 4\mu \parallel \frac{n'}{n^2}\parallel _{\infty }$. And so, $F$ maps $K$ in $K$.   Moreover, $F$ is continuous for the topology of $C([-1,1])$, by definition of $\beta $ in terms of $\alpha $. Finally, $F$ is compact for the topology of $C([-1,1])$, by the compact embedding of $W^{1,\infty }([-1,1])$ in $C([-1,1])$ and the boundedness of $K$ in $W^{1,\infty }([-1,1])$. It follows from a Schauder fixed point theorem that there is a fixed point for $F$ in $K$.\\
Moreover, the solution of (2.5) is unique in the class of non-decreasing functions $\alpha $ such that $\alpha (1)-\alpha (-1)\in [0,\frac{2\mu }{n(1)}]$. Indeed, for any solution $\alpha $ of \eqref{eq:alpha} in this class, $x:= \alpha (1)-\alpha (-1)$ solves $G(x)= 0$, where
\begin{eqnarray*}
G(x):= x-\frac{2}{n(1)}\int _{2\sqrt{\frac{2\mu }{n(1)}}}^{+\infty }\frac{uf_-(u)}{\sqrt{u^2-x}}du+\frac{2}{n(-1)}\int _{-\infty }^{0}\frac{\mid u\mid f_+(u)}{\sqrt{u^2+x}}du\\
-\frac{2}{n(1)}\int _{-\infty }^0f_+(u)du+\frac{2}{n(-1)}\int _{0}^{+\infty }f_-(u)du. 
\end{eqnarray*} 

Using estimates very similar to \eqref{eq:pourG} and $(H_4)$, we can show the function $G$ is increasing. Moreover it follows from $H_3$ and $H_4$ that $G(0)\leq 0$ and $G(\frac{2\mu }{n(1)})\geq 0$. Hence the value of $\alpha (1)-\alpha (-1)$ is unique, as well as $(\alpha (1), \alpha (-1))$, given by
\begin{eqnarray*}
\alpha (1)= \frac{2}{n(1)}\Big( \int _{-\infty }^0f_+(u)du+\int _{2\sqrt{\frac{2\mu }{n(1)}}}^{+\infty }\frac{uf_-(u)}{\sqrt{u^2-\alpha (1)+\alpha (-1)}}du\Big) ,\\
\alpha (-1)= \frac{2}{n(-1)}\Big( \int _{-\infty }^0\frac{\mid u\mid f_+(u)}{\sqrt{u^2+\alpha (1)-\alpha (-1)}}du+\int _0^{+\infty }f_-(u)du\Big) .
\end{eqnarray*}
If $\alpha $ and $\beta $ are two non-decreasing solutions of \eqref{eq:alpha}, then 
\begin{eqnarray*}
(\alpha -\beta )(z)T(z)= 0,\quad z\in [-1,1],
\end{eqnarray*}
where
\begin{eqnarray*}
{ \textstyle
T(z)= \frac{n(z)}{2}\hspace*{1.in} }\\
{ \textstyle -\int _{-\infty }^0\frac{\mid u\mid f_+(u)}{\sqrt{u^2-\alpha (z)+\alpha (1)}\sqrt{u^2-\beta (z)+\alpha (1)}(\sqrt{u^2-\alpha (z)+\alpha (1)}+\sqrt{u^2-\beta (z)+\alpha (1)})}du  }\\
{ \textstyle -\int _{2\sqrt{\frac{2\mu }{n(1)}}}^{+\infty }\frac{uf_-(u)}{\sqrt{u^2-\alpha (z)+\alpha (-1)}\sqrt{u^2-\beta (z)+\alpha (-1)}(\sqrt{u^2-\alpha (z)+\alpha (-1)}+\sqrt{u^2-\beta (z)+\alpha (-1)})}du. }
\end{eqnarray*}
By $H_1$, $H_3$ and $H_4$,
\begin{eqnarray*}
2T(z)\geq n(1)-\int _{-\infty }^0\frac{f_+(u)}{u^2}du-\int _{2\sqrt{\frac{2\mu }{n(1)}}}^{+\infty }\frac{uf_-(u)}{(u^2-\frac{2\mu }{n(1)})^{\frac{3}{2}}}du\\
\geq n(1)-\int _{-\infty }^0\frac{f_+(u)}{u^2}du-\frac{n(1)}{3\sqrt{3}\mu }\int _0^{+\infty }f_-(u)du>0,\quad z\in [-1,1].
\end{eqnarray*}
Hence $\alpha =\beta $.
\begin{theorem}
Assume $H_i$, $1\leq i\leq 4$, and $f_{\pm }\in L^{\infty }$.\\
There is a unique solution $f\in L^{\infty }([-1,1]\times I\! \! R)$ to \eqref{eq:Vlas_red}-\eqref{bc:Vlas_red}, such that $\frac{\rho }{n}\in K$. 
\end{theorem}
\underline{Proof of Theorem 2.3}\hspace*{0.05in}Let $\alpha \in K$ be the solution to \eqref{eq:alpha}. The distribution function $f$ defined by
\begin{eqnarray*}
f(z,v)= f_-(\sqrt{\alpha (z)-\alpha (-1)+v^2}),\quad v>-\sqrt{\alpha (1)-\alpha (z)},\\
f(z,v)= f_+(-\sqrt{\alpha (z)-\alpha (1)+v^2}),\quad v<-\sqrt{\alpha (1)-\alpha (z)},
\end{eqnarray*}
is the unique solution to \eqref{eq:Vlas_red}-\eqref{bc:Vlas_red} in $L^{\infty }([-1,1]\times I\! \! R)$ such that $\frac{\rho }{n}\in K$.

\appendix

\section{Appendix: The polarization operator as \\
a convolution }
\label{App:A}

In this section we restate and prove the result announced in Remark \ref{rem:convol}.

We define
\begin{equation} \label{eq:F0}
F^0_T = \frac{2}{T}   \int_0^{+\infty} \rho e^{-\rho^2/T} (J^0_{\rho})^2 \,d\rho \, ,
\end{equation}
where $T$ is the temperature (of the ions) and $J^0_\rho$ is the gyro-average operator defined in \eqref{eq:J0}. We forget the subscript $L$ in the Larmor radius for conveniance.
We use here the measure $(2\rho/T) e^{-\rho^2/T}$ and not an usual Maxwellain, because we start from a 2D Maxwellian and perform an integration over the angular variable in polar coordinates. 

The operator $F_0$ appears in the electroneutrality equation \eqref{eq:elneut}. Here we will prove the following proposition:
\begin{proposition}
 The operator $F^0_T$ is the convolution with the radial  function $H_T$, defined by
\begin{equation}
H_T(r) = \frac{e^{-\frac{r^2}{4T}}}{2 \pi^{3/2} r \sqrt{T}} \, .
\end{equation}
\end{proposition}

\medskip

\noindent {\bf Proof of the proposition.} 
First, the square of the operator $J_\rho^0$ is
\begin{equation*}
(J^0_\rho)^2(f) (x_g) = \frac{1}{4\pi^2} \int_0^{2\pi} \! \int_0^{2\pi} f(x_g + \rho e^{i\varphi_c} + \rho e^{i\varphi_c'}) \, d\varphi_c d\varphi_c'\, .
\end{equation*}
To simplifiy it, we use the equality 
\begin{equation*}
\rho (e^{i\varphi_c} + e^{i\varphi_c'}) = 2 \rho \cos\left(\frac{\varphi_c-\varphi_c'}{2}\right)e^{i\frac{\varphi_c+\varphi_c'}{2}} \, ,
\end{equation*}
which helps defining the polar coordinates $(r,\phi)$
\begin{equation} \left\{\begin{array}{l}
r=  2 \rho  \left| \cos\left(\frac{\varphi_c-\varphi_c'}{2}\right) \right|, \\
\phi = \frac{\varphi_c+\varphi_c'}{2}  + \ep \, ,
\end{array} \right.   \end{equation}
\bibliographystyle{alpha}
where $\ep = 0$ or $\pi$ depending on the sign of the cosine. This is not exactly a change of variable since a couple $(r,\varphi)$ has exactly two pre-images $(\varphi_c,\varphi_c')$ and  $(\varphi_c',\varphi_c)$. But since it is always bi-valued, we can use the formula for thechange of variables with a factor $2$. The presence of $\ep$ does not introduce specific difficulties, if  the intervals of integration are distinguished. The Jacobian of this transformation is
\begin{equation*}
 \rho \left| \sin\left(\frac{\varphi_c-\varphi_c'}{2}\right) \right|= \sqrt{\rho^2 - \frac{r^2}{4}} \, , 
\end{equation*}
so that
\begin{equation*}
(J^0_\rho)^2(f) (x_g)  = \frac{1}{\pi^2} \int_0^{2\rho} \int_0^{2\pi} f(x_g + r e^{i\phi}) \, \frac{d\phi dr}{\sqrt{4\rho^2 - r^2}} \, , 
\end{equation*}
The next step is to introduce cartesian coordinates. This is a standard transform with Jacobian $r$, so that
\begin{equation*}
(J^0_\rho)^2 (f) (x_g) = \int_{B(2\rho)} f(x_g + y ) \, \frac{dy}{\pi^2 r \sqrt{4\rho^2 - r^2}} \, ,
\end{equation*}
with  $y \in \R^2$, $r=|y|$ and the notation $B(a)$ for the ball of center $0$ and radius $a$. This is exactly the convolution with the radial function
\begin{equation*}
h_\rho(r)  = \frac{\chi_{_{(0,2\rho)}}(r)}{\pi^2 r \sqrt{4\rho^2 - r^2}} \, ,
\end{equation*}
where $\chi_{_A}$ denote the characteristic function of $A$. It can be checked that 
\[
\int h_\rho(r) 2\pi r \, dr =1 \, .
\]

\medskip

The last step is to perform the integration in $\rho$. As $(J_\rho^0)^2$ is the convolution with the radial function $h_\rho$, $ F^0_T = \frac{2}{T} \int \rho e^{-\rho^2/T} (J^0_\rho)^2 \,d\rho$ is the convolution with the radial function
\begin{equation*}
H_T (r) = \frac{2}{ T} \int_0^{+\infty} h_\rho(r) \rho  e^{-\rho^2/T}  \,d\rho \, .
\end{equation*}
In other words,
\begin{equation*}
H_T (r)  =  \frac{2}{\pi^2 r T} \int_{r/2}^{+\infty} \frac{\rho e^{-\rho^2/T}}{ \sqrt{4\rho^2 - r^2}} \,d\rho  \, .
\end{equation*}
We perform the change of variable $ \rho' = \sqrt{4\rho^2 -r^2}$. Hence
\begin{eqnarray*}
H_T (r) & =  & \frac{e^{-\frac{r^2}{4T}}}{2 \pi^2 r T} \int_{0}^{+\infty} e^{-\frac{\rho^2}{4T}} \,d\rho  \\
& = & \frac{e^{-\frac{r^2}{4T}}}{2 \pi^{3/2} r \sqrt{T}} \, . 
\end{eqnarray*}

\medskip

It can be checked that $H_T$ has total mass one. Indeed,
\begin{equation*}
\int_0^\infty  H_T (r) 2\pi r \, dr  =  \frac{1}{ \sqrt{\pi T} } \int_0^\infty  e^{-\frac{r^2}{4T}} \, dr =1
\end{equation*}
\cqfd

{\bibliography{Gyro-refs}} 

\end{document}